\documentclass{amsart}
\usepackage{amsmath,amssymb,verbatim}
\newcommand{\clm} {\mathcal M}
\newtheorem{thm}{Theorem}

\theoremstyle{definition}

\newtheorem{remark}{Remark}

\newcommand{\e}{\epsilon}
\newcommand{\no}{\noindent}
\newcommand{\nah}{\nonumber}

\newcommand{\hatpsi}{\widehat\psi}
\newcommand{\hatphi}{\widehat\phi}
\newcommand{\la}{\lambda}
\newcommand{\ds}{\displaystyle}

\begin{document}
\title[Roots of Toeplitz Operators]{Roots of Toeplitz Operators on the Bergman space}

\date{\today}

\author{I. Louhichi \and N. V. Rao}
\address{Department of Mathematics and Statistics\\
King Fahd University of Petroleum and Minerals\\
Dhahran, Saudi Arabia}
\address{The University Of Toledo, College of Arts and Sciences,
Department of Mathematics, Mail Stop 942. Toledo, Ohio 43606-3390,
USA}
\email{issam@kfupm.edu.sa}
\email{rnagise@math.utoledo.edu}
\begin{abstract}
One of the major questions in the theory of Toeplitz operators on the Bergman space over the
unit disk $\mathbb D$ in the complex plane $\mathbb C$ is a complete description of the commutant
of a given Toeplitz operator, that is
the set of all Toeplitz operators that commute with it. In \cite{l}, the first author obtained a
complete description
of the commutant of Toeplitz operator $T$ with any quasihomogeneous symbol $\phi(r)e^{ip\theta}, p>0$
in case it has a Toeplitz p-th root $S$ with symbol $\psi(r)e^{i\theta}$, namely,
commutant of $T$ is the closure of the linear space generated by powers $S^n$ which are Toeplitz.
But the existence of p-th root was known until now
only when $\phi(r)=r^m,m \geq 0$. In this paper we will show the existence of p-th roots for a much
larger class of symbols, for example, it includes such symbols for which
$$\phi(r)=\sum_{i=1}^kr^{a_i}(\ln r)^{b_i},0\leq a_i, b_i\mbox{ for all }  1\leq i\leq k .$$

\end{abstract}

\maketitle
\section{Introduction}
\noindent Let $\mathbb{D}$ be the unit disc in the complex plane
$\mathbb{C}$, and $dA=rdr\frac{d\theta}{\pi}$ be the normalized
Lebesgue area measure so that the measure of $\mathbb{D}$ equals
1. Let  $L^2_a$ be The Bergman space, the Hilbert space of
 functions, analytic on $\mathbb{D}$ and square
integrable with respect to the measure $dA$. We denote the inner
product in $L^2(\mathbb{D},dA)$ by $<,>$. It is well known that
$L^2_a$ is a closed subspace of the Hilbert space
$L^2(\mathbb{D},dA)$, with the set of functions
$\{\sqrt{n+1}z^n\,\mid\,n\geq 0\}$ as an orthonormal basis. Let
$P$ be the orthogonal projection from $L^2(\mathbb{D},dA)$ onto
$L^2_a$. For a bounded function $f$ on $\mathbb{D}$, the
Toeplitz operator $T_{f}$ with symbol $f$ is defined by
$$T_{f}(h)=P(fh) \mbox { for } h\in L^2_a.$$

A symbol $f$ is said to be quasihomogeneous of order $p$ an integer, if it can be written as $f(re^{i\theta})=e^{ip\theta}\phi(r)$, where $\phi$ is a radial function on $\mathbb{D}$. In this case, the associated Toeplitz operator $T_f$ is also called quasihomogeneous Toeplitz of order $p$. Quasihomogeneous Toeplitz operators were first introduced in \cite{lz} while generalizing the results of \cite{cr}. We assume $p>0$ from now on.

We are looking for, given a quasihomogeneous operator $T$ of degree $p$, a quasihomogeneous operator $S$ of degree 1 such that $S^p=T$. It was proved in \cite{l}
that any such root if it exists, is unique up to a multiplicative constant. Also the existence of $p-$th roots for the case $\phi(r)=r^m$ for any arbitrary $m\geq 0,p>0$ was proved in \cite{l} using the results in \cite{cr}. Here we plan to deal with more general $\phi(r)$.
\section{ The Mellin Transform and two Lemmas}
For any two functions $f(r)$ and $g(r)$ defined on $I=[0,1]$, we define the Mellin convolution  as follows:
$$(f*g)(r)=\int_r^1f(\frac{r}{t})g(t)\frac{dt}{t}.$$ Often we are interested in knowing when the convolution is
a bounded function in the interval $I$. To that purpose, we introduce the following concept of the
type for a function f. We say f is of type $(a,b)$ with $a\geq 0$ and $b>0$ if
$$|f(r)|\leq Cr^a(1-r)^{b-1}$$ on $I$, where $C$ is a constant depending on $f$. Also we express the same thing as $$f(r)\ll r^a(1-r)^{b-1}$$ omitting the constants and the absolute value signs.

\no{\bf Lemma A}. Suppose $f(r)$ is of type $(a,b)$ and $g(r)$ is of type $(c,d)$. Then  their convolution product
\begin{displaymath}\left\{\begin{array}{ll}
(f*g)\textrm{ is of type } \left(\min\{a,c\},b+d\right)& \textrm{ if } a\neq c\\
\textrm{ and }&\\
(f*g)(r)\ll r^{\min\{a,c\}}(1-r)^{b+d-1}\ln(\frac{e}{r})&\textrm{ if } a=c.\end{array}\right.
\end{displaymath}
 This can be generalized to any finite product as follows: Suppose for $1\leq i\leq n$, $f_i(r)$ is of type $(a_i,b_i)$.
Then $h(r)$, their convolution product satisfies
$$h(r)\ll r^{\alpha}(1-r)^{\beta-1}\left(\ln\left(\frac er\right)\right)^{n-1}\eqno(1)$$
where $\alpha=\min\{a_i\},\beta=\sum b_i$. Further, if we know the number of $a_i$ that are equal to $\min\{a_i\}$ to be say $l$, the estimate (1) can be improved to $$h(r)\ll r^{\alpha}(1-r)^{\beta-1}\left(\ln\left(\frac er\right)\right)^{l-1}.\eqno(2)$$ Thus the log term will disappear if $l=1$.

\begin{remark}  Most of the time our aim is to prove $h$ is bounded and the presence of log does not interfere with that aim since $\alpha>0$ which makes $h(r)$ bounded near zero and since $\beta\geq 1$, it is bounded near 1. But log cannot be avoided. Take for example $f_i(r)=r$ for every $i$ and compute the
convolution product. It checks out to be $\ds \frac{r(\ln r)^{n-1}}{(n-1)!}$, by a simple integration.
\end{remark}
\no {\bf Lemma B}. Suppose $\ds f_i(r)=r^{a_i}(1-r)^{b_i-1}$ where $a_i,b_i$ are positive for $1\leq i\leq n$.
Let $\alpha,\beta$ be as defined in Lemma A. Given any integer $k\geq 0$, the $k-$the derivative of $h$,
the convolution product of $f_i$,
satisfies the following: $$h^{(k)}(r)\ll r^{\alpha-k}(1-r)^{\beta-k-1}\left(\ln\left(\frac{e}{r}\right)\right)^{n-1}.$$ Here the constant  involved depends on $k$ and $h$.
\section{applications of Lemmas A and B}
One of our most useful tools in the following calculations will be the Mellin transform.
The Mellin transform $\widehat{\phi}$
of a radial function $\phi$ in $L^1([0,1],rdr)$ is defined by
$$\widehat{\phi}(z)=\int_{0}^{1} f(r) r^{z-1}\,dr=\mathcal{M}(\phi)(z).$$
It is well known that, for these functions, the Mellin transform
is well defined on the right half-plane $\{z : \Re z\geq 2\}$ and
it is analytic on $\{z:\Re z>2\}$. It is important and helpful to
know that the Mellin transform $\widehat{\phi}$ is uniquely
determined by its values on any arithmetic sequence of integers.
In fact we have the following classical theorem \cite[p.102]{rem}.

\begin{thm}
\label{blk1}
 Suppose that $f$ is a bounded analytic function on
$\{ z : \Re z>0\}$ which vanishes at the pairwise distinct points
$z_1, z_2 \cdots$, where
\begin{itemize}
\item[i)] $\inf\{|z_n|\}>0$\\
 and
\item[ii)] $\sum_{n\ge 1}\Re(\frac{1}{z_n})=\infty$.
\end{itemize}
Then $f$ vanishes identically on $\{ z : \Re z>0\}$.
\end{thm}
\begin{remark}
\label{blk2} Now one can apply this theorem to prove that if $\phi\in
L^1([0,1],rdr)$ and if there exist $n_0, p\in\mathbb{N}$ such that
$$\widehat{\phi}(pk+n_0)=0\textrm{ for all } k\in\mathbb{N},$$
then $\widehat{\phi}(z)=0$ for all $z\in\{z:\Re z>2\}$ and so $\phi=0$.
\end{remark}
Moreover, it is easy to see that the Mellin transform converts the convolution product into a pointwise product, i.e that:
$$\widehat{(\phi\ast \psi)}(r)=\widehat{\phi}(r)\widehat{\psi}(r).$$
A direct calculation shows that a quasihomogeneous Toeplitz operator acts on the elements of the orthogonal basis of $L^2_a$ as a shift operator with a holomorphic weight. In fact, for $p\geq 0$ and for all $k\geq 0$, we have
\begin{eqnarray*}
T_{e^{ip\theta}\phi}(z^k)&=&P(e^{ip\theta}\phi z^k)=\sum_{n\geq
0}(n+1)\langle e^{ip\theta}\phi z^k, z^n\rangle z^n\\
&=& \sum_{n\geq 0}(n+1)
\int_{0}^{1}\int_{0}^{2\pi}\phi(r)r^{k+n+1}e^{i(k+p-n)\theta}
\frac{d\theta}{\pi}dr z^n\\
&=&2(k+p+1)\widehat{\phi}(2k+p+2)z^{k+p}.
\end{eqnarray*}
Now we are ready to start with the following relatively easy example.

\subsection{$p$-th roots of $T_{e^{ip\theta}\phi}$ where $\phi(r)=r+r^2$.}

Does there exist a radial function $\psi$ such that $\left(T_{e^{i\theta}\psi}\right)^p=T_{e^{ip\theta}\phi}$? If it is the case, then we will have
$$\left(T_{e^{i\theta}\psi}\right)^p(z^k)=T_{e^{ip\theta}\phi}(z^k), \textrm{ for all } k\geq0.$$
Since
$$\left(T_{e^{i\theta}\psi}\right)^p(z^k)=
\left[\prod_{j=0}^{p-1}(2k+2j+4)\widehat{\psi}(2k+2j+3)\right]z^{k+p},$$
we obtain for all integers $k\geq 0$
$$(2k+2p+2)\widehat{\phi}(2k+p+2)=
\left[\prod_{j=0}^{p-1}(2k+2j+4)\widehat{\psi}(2k+2j+3)\right],$$ from which and Remark \ref{blk2}
follows,, by setting $z=2k+3$, the identity, valid in the right halfplane

\begin{equation}\label{xxx}
(z+2p-1)\widehat{\phi}(z+p-1)=
\left[\prod_{j=0}^{p-1}(z+2j+1)\widehat{\psi}(z+2j)\right].
\end{equation}

If we divide the equation (\ref{xxx}) by the equation obtained by replacing $z$ by $z+2$ in the equation (\ref{xxx}), after cancelation, we obtain that in the right halfplane,

\begin{equation}\label{main}
\frac{\widehat{\psi}(z+2p)}{\widehat{\psi}(z)}=
\frac{(z+1)\widehat{\phi}(z+p+1)}{(z+2p-1)\widehat{\phi}(z+p-1)}, \textrm{ for } \Re z>0.
\end{equation}
Since $\displaystyle\widehat{\phi}(z)=\frac{1}{z+1}+\frac{1}{z+2}=\frac{2z+3}{(z+1)(z+2)}$, it follows that
$$\frac{\widehat{\psi}(z+2p)}{\widehat{\psi}(z)}=\frac{(z+1)}{(z+2p-1)}\frac{(2z+2p+5)}{(z+p+2)(z+p+3)}\frac{(z+p)(z+p+1)}
{(2z+2p+1)}, \textrm{ for }\Re z>0.$$
If we denote by $\la(\zeta)=\hatpsi(2p\zeta)$, the above equation becomes
$$\frac{\la(\zeta +1)}{\la(\zeta)}=\frac{(2p\zeta+1)(4p\zeta+2p+5)(2p\zeta+p)(2p\zeta+p+1)}
{(2p\zeta+2p-1)(2p\zeta+p+2)(2p\zeta+p+3)(4p\zeta+2p+1)}, \textrm{ for } \Re \zeta>0.$$
Using the well-known identity $\Gamma(z+1)=z\Gamma(z)$, where $\Gamma$ is the Gamma function, we can write that
\begin{equation}\label{gam}
\frac{\la(\zeta +1)}{\la(\zeta)}=\frac{F(\zeta+1)}{F(\zeta)}\textrm{ for } \Re \zeta>0,
\end{equation}
where
$$F(\zeta)=\frac{\Gamma(\zeta+a_1)\Gamma(\zeta+a_2)\Gamma(\zeta+a_3)\Gamma(\zeta+a_4)}
{\Gamma(\zeta+a'_1)\Gamma(\zeta+a_2')\Gamma(\zeta+a_3')\Gamma(\zeta+a_4')},$$
with $a_i$ are in increasing order $\displaystyle{\frac2{4p}, \frac{2p}{4p}, \frac{2p+2}{4p}, \frac{2p+5}{4p}}$ respectively and $a_i'$ are in almost increasing order
$\displaystyle{\frac{2p+1}{4p}, \frac{2p+4}{4p}, \frac{4p-2}{4p}, \frac{2p+6}{4p}}$ respectively for $i=1,\ldots ,4$. Equation (\ref{gam}), combined with \cite[Lemma 6,~p.1428]{l}, gives us that there exists a constant $C$ such that
\begin{equation}\label{equal}
\la(\zeta)=CF(\zeta), \textrm{ for }\Re\zeta>0.
\end{equation}
Basic observation is that the quotient of two Gamma functions
$$\frac{\Gamma(\zeta+a_i)}{\Gamma(\zeta+a_i')}, \textrm{ where }0<a_i<a_i'$$ is a constant times the Beta function
$$B(\zeta+a_i,a_i'-a_i)=\int_0^1x^{\zeta+a_i-1}(1-x)^{a_i'-a_i-1}\,dx.$$ Moreover, according to our definition of the Mellin transform, it turns out that $B(\zeta+a_i,a_i'-a_i)$ is the Mellin Transform of
$x^{a_i}(1-x)^{a_i'-a_i-1}$ which is of type $(a_i,a_i'-a_i)$. Since the $a_i$ are smaller than $a'_i$
respectively for $i=1,\ldots ,4$ (in fact $a_3'\geq a_3$ if and only if $2p\geq 4$ which is always true), Equation (\ref{equal}) implies that
$$\lambda(\zeta)=C\Pi_{i=1}^{4}B(\zeta+a_i,a'_i-a_i),$$
where $C$ is a constant. Since the product of Mellin transforms equals to the Mellin of the convolution product, we would have
$$\lambda(\zeta)=Ch(\zeta),$$
where $h$ is the convolution product of four functions of type $(a_i,a_i'-a_i)$, $i=1,\ldots ,4$. Now Lemma A tells us that
$$h(r)\ll r^{\min\{a_i\}}(1-r)^{\sum_i(a_i'-a_i) -1}\ln(\frac{e}{r}).$$
 Because $\sum_ia'_i-a_i=1$, we have
 $$h(r)\ll r^{\min\{a_i\}}\ln(\frac{e}{r}),$$ and hence $h$ is bounded function.
 Therefore the function $\psi$, if it exists,
 satisfies the equation $$\hatpsi(2p\zeta)=C\widehat h(\zeta)$$ for some constant $C$,  which is equivalent to
 $$\int_0^1\psi(r)r^{2p\zeta-1}dr=C\int_0^1h(t)t^{\zeta-1}dt.$$
 Now, by a change of variables $t=r^{2p}$, we obtain
 $$\int_0^1\psi(r)r^{2p\zeta-1}dr=\int_0^1h(r^{2p})r^{2p\zeta-1}2pdr.$$
Thus $\ds\psi(r)=2ph(r^{2p})$, and so $\psi$ is bounded. Hence the operator $T_{e^{i\theta}\psi}$ is a genuine Toeplitz operator and $p-$th root of $T_{e^{ip\theta}\phi}$.

\subsection{$p$-th roots of $T_{e^{ip\theta}\phi}$ where $\widehat{\phi}(z)$ is a proper rational fraction.}
We recall that if there exists a radial function $\psi$ such that $\left(T_{e^{i\theta}\psi}\right)^p=T_{e^{ip\theta}\phi}$, then we have Equation (\ref{main}) which is
$$\hatpsi(z+2p)=\hatpsi(z)\frac{(z+1)\hatphi(z+p+1)}{(z+2p-1)\hatphi(z+p-1)}, \textrm{ for } \Re z>0. $$
Here we are assuming $\widehat{\phi}(z)=\frac{P(z)}{Q(z)}$ where $P(z)=\ds\prod_{j=1}^m(z+a_j)$ and $Q(z)=\ds\prod_{k=1}^n(z+b_k)$ with $1\leq m<n$. So that
\begin{eqnarray}\hatpsi(z+2p)&=&\hatpsi(z)\frac{(z+1)}{(z+2p-1)}\frac{P(z+p+1)Q(z+p-1)}{P(z+p-1)Q(z+p+1)}\nah\\
&=&\frac{(z+1)}{(z+2p-1)} \prod_{j=1}^m\frac{z+a_j+p+1}{z+a_j+p-1}\prod_{k=1}^n\frac{z+b_k+p-1}{z+b_k+p+1}\nah
\end{eqnarray}
Let $\la(\zeta)=\hatpsi(2p\zeta)$. Then the equality above becomes
$$\frac{\la(\zeta +1)}{\la(\zeta)}=\frac{(2p\zeta +1)}{(2p\zeta +2p-1)} \prod_{j=1}^m\frac{2p\zeta +a_j+p+1}{2p\zeta +a_j+p-1}\prod_{k=1}^n\frac{2p\zeta +b_k+p-1}{2p\zeta +b_k+p+1}$$
Therefore, by \cite[Lemma 6,~p.1428]{l}, $\lambda$ is constant times the quotient of $m+n+1$ Gamma functions in the numerator and about the same in the denominator as follows:
\begin{equation}\label{multigam}
\la(\zeta)=C\frac{\Gamma(\zeta +A_0)}{\Gamma(\zeta +A'_0)} \prod_{j=1}^m\frac{\Gamma(\zeta +A_j)}{\Gamma(\zeta +A'_j)}\prod_{k=1}^n\frac{\Gamma(\zeta +B_k)}{\Gamma(\zeta +B'_k)}
\end{equation}
where $A_0=\frac1{2p}$, $A'_0=\frac{2p-1}{2p}$, $A_j=\frac{a_j+p+1}{2p}$, $A'_j=\frac{a_j+p-1}{2p}$,
 $B_k=\frac{b_k+p-1}{2p}$ and $B'_k=\frac{b_k+p+1}{2p}$ for $1\leq j\leq m$ and $1\leq k\leq n$.
 Based on the same argument as in the previous subsection, we would like to write each quotient of two Gamma functions as a constant times a Beta function. In order to do that, we must assume that all $A_j$ and $B_k $ are positive for every $0\leq j\leq m$ and $1\leq k\leq n$. Moreover we observe that
 $$A'_0-A_0=\ds\frac{p-1}p,\  A'_j-A_j=-\frac1p,\  B'_k-B_k=\frac1p.$$
 So each quotient of two Gamma functions in Equation (\ref{multigam})can be written as a constant times a Beta function except those involving $A_j$ for $1\leq j\leq m$. We fix this matter by noting that
$\Gamma(\zeta +A'_j+1)=(\zeta +A'_j)\Gamma(\zeta +A'_j)$, and so here $A'_j+1-A_j=\frac{p-1}p$. Hence, Equation (\ref{multigam}) becomes
$$\frac{\lambda(\zeta)}{\prod_{j=1}^m(\zeta+A'_j)}=C\frac{\Gamma(\zeta +A_0)}{\Gamma(\zeta +A'_0)} \prod_{j=1}^m\frac{\Gamma(\zeta +A_j)}{\Gamma(\zeta +A'_j+1)}\prod_{j=1}^n\frac{\Gamma(\zeta +B_j)}{\Gamma(\zeta +B'_j)}.$$
As in the previous subsection, this quotient of $m+n+1$ Gamma functions on the numerator and the same in the denominator,
respectively would be the Mellin transform of the convolution product of $m+n+1$ functions. Let us denoted it $h$. By Lemma A, we have
$$h(r)\ll r^A(1-r)^{B-1}\left(\ln\left(\frac{e}{r}\right)\right)^{m+n},$$ where $A=\min\{A_j\}$ which is definitely positive, and $B$ is given by $$\ds A'_0-A_0+\sum_{j=1}^m A'_j+1-A_j+\sum_{k=1}^n B'_k-B_k=(m+1)\frac{p-1}p+\frac np=m+1+\frac{n-m-1}p.$$
 Therefore we obtain
 \begin{equation*}
 h(r)\ll r^A(1-r)^{m+\frac{n-m-1}p}\left(\ln\left(\frac er\right)\right)^{m+n}=r^A(1-r)^{m+\upsilon}\left(\ln\left(\frac er\right)\right)^{m+n},
 \end{equation*}
 where $\upsilon=\ds\frac{n-m-1}p$ is a non-negative number. Using Lemma B, we see that $h$ has all derivatives of order not exceeding $m$ and they satisfy the inequality
 \begin{equation*}
 r^jh^{(j)}(r)\ll r^A(1-r)^{m-j+\upsilon}\left(\ln\left(\frac er\right)\right)^{m+n}.
 \end{equation*}
 Further  the function $\psi$, if it exists, would
 satisfy the equation
 \begin{equation} \label{psi}
 \hatpsi(2p\zeta)=C\left(\prod_{j=1}^m(\zeta+A'_j)\right)\widehat h(\zeta).
 \end{equation}
 Now it is easy to check by integration by parts the following identity
 \begin{equation*}
 \zeta\widehat h(\zeta)= -\clm\left(r\frac{dh}{dr}\right)(\zeta)
 \end{equation*}
 provided $h$ vanishes at 1 and $rh'$ is bounded in $(0,1)$. Thus in the current case, denoting  $h'$ by $Dh$ where $D=\ds\frac d{dr}$, we can see
 $$(\zeta+A'_j)\widehat h(\zeta)= \clm\left(\left(A'_j-rD\right)h\right)(\zeta), $$ and
 $$\left(\prod_{j=1}^m(\zeta+A'_j)\right)\widehat h(\zeta) = \clm\left(\prod_{j=1}^m(A'_j-rD)h\right)(\zeta).$$ Let us set
 $$H(r)=\left(\prod_{j=1}^m(A'_j-rD)h\right)(r)$$
 which allows us to rewrite Equation (\ref{psi}) as
 $$\int_0^1\psi(r)r^{2p\zeta-1}dr=C\int_0^1H(t)t^{\zeta-1}\,dt.$$
 Now, by a change of variables $t=r^{2p}$, we obtain
 $$\int_0^1\psi(r)r^{2p\zeta-1}dr=C\int_0^1H(r^{2p})r^{2p\zeta-1}2pdr.$$
 Thus $\ds\psi(r)=2pCH(r^{2p})$, and hence is bounded and the operator $T_{e^{i\theta}\psi}$ is a genuine Toeplitz operator and $p-$th root of $T_{e^{ip\theta}\phi}$.

\section{Proof of the Lemma A for two functions}

We choose to start proving Lemma A for two functions $f$ and $g$ of type $(a,b)$ and $(c,d)$ respectively, with $a$, $b$, $c$ and $d$ are all positive. Similar thing was discussed in \cite[pages 210-212]{cr} but with less generality since the goal was different.

Let $h(r)=(f*g)(r)$. By definition of the Mellin convolution, it is easy to see that
$$h(r)\ll\int_r^1\left(\frac rt\right)^a\left(1-\frac rt\right)^{b-1}t^c(1-t)^{d-1}\,\frac{dt}{t},$$ which after a change of variables
$\ds\frac{t-r}{1-r}=u$ and using the consequent identities
$$t=r+u-ru=,\  t-r=u(1-r),\ 1-t=(1-u)(1-r),\  dt=(1-r)du$$
leads to
\begin{eqnarray*}
h(r)&\ll&\int_r^1\left(\frac rt\right)^a\left(1-\frac rt\right)^{b-1}t^c(1-t)^{d-1}\,\frac{dt}{t}\\
&=&\int_r^1\left(\frac rt\right)^a\left(\frac{t-r}{t}\right)^{b-1}t^c(1-t)^{d-1}\,\frac{dt}{t}\\
&=&\int_0^1r^at^{-a}u^{b-1}(1-r)^{b-1}t^{-b+1}t^c(1-u)^{d-1}(1-r)^{d-1}(1-r)\frac{du}{t}\\
&=& r^a(1-r)^{b+d-1}\int_0^1t^{c-a-b}u^{b-1}(1-u)^{d-1}du.
\end{eqnarray*}
We have the following cases
\begin{itemize}
\item[$\bullet$] If $c-a-b\geq 0$. Since $0\leq t\leq 1$, we have
$$h(r)\ll r^a(1-r)^{b+d-1},$$
and hence $h$ is of type $(a,b+d)$.
\item[$\bullet$] If $c-a-b<0$. Assuming $c-a>0$ and noting that $t\geq u$, we obtain
\begin{eqnarray*}
h(r)&\ll& r^a(1-r)^{b+d-1}\int_0^1u^{c-a-b}u^{b-1}(1-u)^{d-1}\,du\\
&\leq& r^a(1-r)^{b+d-1}\int_0^1u^{c-a-1}(1-u)^{d-1}\,du\\
&=&  r^a(1-r)^{b+d-1}B(c-a,d),
\end{eqnarray*}
and therefore $h$ is of type $(a,b+d)$.

Now in case $c=a$, for any  number $0<\e\leq b$,  noticing that $t\geq r$ and  $u>0$, we have
\begin{eqnarray*}
h(r)&=& r^a(1-r)^{b+d-1}\int_0^1t^{-b}u^{b-1}(1-u)^{d-1}\,du\\
&\ll& r^a(1-r)^{b+d-1}\int_0^1t^{-\e}t^{\e-b}u^{b-1}(1-u)^{d-1}\,du\\
&\leq& r^a(1-r)^{b+d-1}\int_0^1r^{-\e}u^{\e-b}u^{b-1}(1-u)^{d-1}\,du\\
&\leq& r^a(1-r)^{b+d-1}r^{-\e}\int_0^1u^{\e-1}(1-u)^{d-1}\,du\\
&\leq& r^a(1-r)^{b+d-1}B(\e,d)r^{-\e}.
\end{eqnarray*}
Now since $\e B(\e,d)= \frac {\Gamma(\e+1)\Gamma(d)}{\Gamma(\e+d)}$ is holomorphic as
a function of $\e$ in a neighborhood of the interval $(0,b)$, there exists a constant $C$ such that
$\e B(\e,d)\leq C$ on that interval, and therefore
$$h(r)\leq C r^a(1-r)^{b+d-1}r^{-\e}\e^{-1},\textrm{ for every } 0<\e\leq b.$$
Here we emphasize the fact that $C$ does not depend on $r$ and $\e$ as long as $0<r<1$ and $0<\e\leq b$. For a fixed but arbitrary $r$, let $E(\epsilon)=r^{-\e}\e^{-1}$ and $m(r)=\ds\min_{(0,b]}E(\e)$.  Then
\begin{equation}\label{m}
h(r)\leq C r^a(1-r)^{b+d-1}m(r).
\end{equation}
Moreover the function $E$ decreases in the interval $\left(0,-\frac{1}{\ln r}\right)$ and increases
in the interval $\left(-\frac{1}{\ln r},+\infty\right)$. Further $-\frac{1}{\ln r}\leq b$ if and only if $r\leq e^{-\frac{1}{b}}$. Thus Equation (\ref{m}) implies
\begin{itemize}
\item[]If $r\leq e^{-\frac{1}{b}}$,
\begin{equation}\label{l}
h(r)\ll r^a(1-r)^{b+d-1}m(r)\leq r^a(1-r)^{b+d-1}e\ln\left(\frac{1}{r}\right).
\end{equation}
\item[]If $r>e^{-\frac{1}{b}}$,
\begin{equation}\label{g}
h(r)\ll r^a(1-r)^{b+d-1}r^{-b}b^{-1}\leq r^a(1-r)^{b+d-1}\frac eb.
\end{equation}
\end{itemize}
Combining (\ref{l}) and (\ref{g}), we obtain
\begin{eqnarray*}
h(r)&\ll& r^a(1-r)^{b+d-1}\left(e\ln\left(\frac{1}{r}\right)+\frac eb\right)\\
&\ll& r^a(1-r)^{b+d-1}\ln\left(\frac{e}{r}\right),\textrm{ for all }0<r<1.
\end{eqnarray*}
\end{itemize}

\section{Lemma A for convolution product of more than two functions}
In this context we can assume that the function $f_i$ which is of type $(a_i,b_i)$ is
$$f_i(x)=x^{a_i}(1-x)^{b_i-1},\textrm{ for } 1\leq i\leq n.$$
The convolution product of these $n$ functions is defined by a repeated integral
\begin{eqnarray}\label{int}
h(r)&=&\int_r^1\int_{r/x_1}^1\int_{r/x_1x_2}^1\ldots\int_{r/x_1x_2\ldots x_{n-2}}^1\nah\\
&&f_1(x_1)f_2(x_2)\ldots f_{n-1}(x_{n-1})f_n\left(\frac r{x_1\ldots x_{n-1}}\right)
\frac{dx_{n-1}}{x_{n-1}}\ldots\frac{dx_3}{x_3}\frac{dx_2}{x_2}\frac{dx_1}{x_1}
\end{eqnarray}
As in the case of two functions where we made a change of variables
 $u=\frac{t-r}{1-r},$ we make change of variables so that the new integral is over the unit cube
$I^{n-1}$ where limits of integration do not depend on other variables.
Let $y_0=1$ and inductively define $y_i=\prod_{j=1}^ix_i$ for $i\geq 1$. Now we make the change of variables as
follows: $$x_i=\frac r{y_{i-1}}+\left(1-\frac r{y_{i-1}}\right)\xi_i,\textrm{ for }i\geq1,$$ so that the limits for each $\xi_i$ are $0$ and $1$. Further we note
\begin{eqnarray*}
y_i-r=x_iy_{i-1}-r=(y_{i-1}-r)\xi_i,\textrm{ for }i\geq0.
\end{eqnarray*}
Let us set $\eta_0=1$ and $\eta_i=\prod_{j=1}^i\xi_i, \textrm{ for }i\geq1.$ It is easy to show, by induction on $i$, that
$$y_i-r=(1-r)\eta_i, \textrm{ for all } i\geq 1.$$
Further
\begin{equation}\label{x}
(1-x_i)=(1-\xi_i)\left(1-\frac r{y_{i-1}}\right)=\frac{(1-\xi_i)(1-r)\eta_{i-1}}{y_{i-1}},\textrm{ for all } i\geq 1.
\end{equation}
Thus
$$f_i(x_i)=x_i^{a_i}(1-x_i)^{b_i-1}=\left(\frac{y_i}{y_{i-1}}\right)^{a_i}\left(\frac{(1-\xi_i)
(1-r)\eta_{i-1}}{y_{i-1}}\right)^{b_i-1}, \textrm{ for } 1\leq i\leq n-1.$$
But for $i=n$, we have
\begin{eqnarray*}
f_n(x_n)&=&\left(\frac r{y_{n-1}}\right)^{a_n}\left(1-\frac r{y_{n-1}}\right)^{b_n-1}\\
   &=&\left(\frac r{y_{n-1}}\right)^{a_n}\left(\frac{y_{n-1}- r}{y_{n-1}}\right)^{b_n-1}\\
 &=&\left(\frac r{y_{n-1}}\right)^{a_n}\left(\frac{(1-r)\eta_{n-1}}{y_{n-1}}\right)^{b_n-1}.
\end{eqnarray*}
Writing the product of functions in (\ref{int}) in terms of $\xi_i$, $\eta_i$, $r$, $y_i$, for $1\leq i \leq n-1$ yields
\begin{equation}\label{prod}
r^{a_n}\prod_{i=1}^{n-1}\eta_i^{b_{i+1}-1}\prod_{i=1}^{n-1}y_i^{a_i-a_{i+1}-b_{i+1}+1}
\prod_{i=1}^{n-1}(1-\xi_i)^{b_i-1}\prod_{i=1}^{n}(1-r)^{b_i-1}.
\end{equation}
Using equalities of $(\ref{x})$, we calculate the differential form
\begin{equation}\label{diff}
\bigwedge_{i=1}^{n-1}\frac{dx_i}{x_i}=\bigwedge_{i=1}^{n-1}\frac{(1-r)\eta_{i-1}d\xi_i}{x_iy_{i-1}}
=(1-r)^{n-1}\prod_{i=1}^{n-1}\frac{\eta_{i-1}}{y_i}\bigwedge_{i=1}^{n-1}d\xi_i.
\end{equation}
From $(\ref{int}),(\ref{x}),$ and $(\ref{prod})$ we derive
\begin{eqnarray}\label{h}
h(r)&=&r^{a_n}(1-r)^{b_1+\ldots+b_n-1}\int_{I^{n-1}}\eta_{n-1}^{b_n-1}
\prod_{i=1}^{n-2}\eta_i^{b_{i+1}}\nah\\
& &\prod_{i=1}^{n-1}y_i^{a_i-a_{i+1}-b_{i+1}}
\prod_{i=1}^{n-1}(1-\xi_i)^{b_i-1}\bigwedge_{i=1}^{n-1}d\xi_i.
\end{eqnarray}
Let us assume that $a_i$ are arranged in decreasing order. Then the product $$\eta_i^{b_{i+1}}y_i^{a_i-a_{i+1}-b_{i+1}}\leq 1$$ since
$\eta_i\leq y_i\leq 1$. Therefore $$h(r)\leq r^{a_n}(1-r)^{b_1+\ldots+b_n-1}\int_{I^{n-1}}\eta_{n-1}^{b_n-1}
y_{n-1}^{a_{n-1}-a_n-b_n}\prod_{i=1}^{n-1}(1-\xi_i)^{b_i-1}\bigwedge_{i=1}^{n-1}d\xi_i.$$
Here four cases have to be discussed.
\begin{itemize}
\item[\bf{Case 1.}] If $a_{n-1}-a_n-b_n\geq 0$. Then
\begin{eqnarray*}
h(r)&\leq &r^{a_n}(1-r)^{b_1+\ldots+b_n-1}\int_{I^{n-1}}\eta_{n-1}^{b_n-1}
\prod_{i=1}^{n-1}(1-\xi_i)^{b_i-1}\bigwedge_{i=1}^{n-1}d\xi_i\\
&\leq& r^{a_n}(1-r)^{b_1+\ldots+b_n-1}\prod_{i=1}^{n-1}\frac{\Gamma(b_n)\Gamma(b_i)}{\Gamma(b_n+b_i)}.
\end{eqnarray*}
\item[\bf{Case 2.}] If $a_{n-1}-a_n-b_n< 0$ and $a_{n-1}\neq a_n$. Then
\begin{eqnarray*}
h(r)&\leq& r^{a_n}(1-r)^{b_1+\ldots+b_n-1}\int_{I^{n-1}}y_{n-1}^{a_{n-1}-a_n-b_n}\eta_{n-1}^{b_n-1}
\prod_{i=1}^{n-1}(1-\xi_i)^{b_i-1}\bigwedge_{i=1}^{n-1}d\xi_i\\
&\leq& r^{a_n}(1-r)^{b_1+\ldots+b_n-1}\int_{I^{n-1}}\eta_{n-1}^{a_{n-1}-a_n-b_n}\eta_{n-1}^{b_n-1}
\prod_{i=1}^{n-1}(1-\xi_i)^{b_i-1}\bigwedge_{i=1}^{n-1}d\xi_i\\
&\leq& r^{a_n}(1-r)^{b_1+\ldots+b_n-1}\int_{I^{n-1}}\eta_{n-1}^{a_{n-1}-a_n-1}
\prod_{i=1}^{n-1}(1-\xi_i)^{b_i-1}\bigwedge_{i=1}^{n-1}d\xi_i\\
&\leq& r^{a_n}(1-r)^{b_1+\ldots+b_n-1}\prod_{i=1}^{n-1}\frac{\Gamma(a_{n-1}-a_n)\Gamma(b_i)}{\Gamma(a_{n-1}-a_n+b_i)}.
\end{eqnarray*}
\item[\bf{Case 3.}] If $a_{n-1}=a_n$. Choose an arbitrary $0<\e\leq b_n$ and note that $y_{n-1}\geq r$, and $\eta_{n-1}>0$. Therefore
\begin{eqnarray}
h(r)&\leq& r^{a_n}(1-r)^{b_1+\ldots+b_n-1}\int_{I^{n-1}}y_{n-1}^{-b_n}\eta_{n-1}^{b_n-1}
\prod_{i=1}^{n-1}(1-\xi_i)^{b_i-1}\bigwedge_{i=1}^{n-1}d\xi_i\nah\\
&\leq& r^{a_n}(1-r)^{b_1+\ldots+b_n-1}\int_{I^{n-1}}y_{n-1}^{-\e}y_{n-1}^{\e-b_n}\eta_{n-1}^{b_n-1}
\prod_{i=1}^{n-1}(1-\xi_i)^{b_i-1}\bigwedge_{i=1}^{n-1}d\xi_i\nah\\
&\leq& r^{a_n}(1-r)^{b_1+\ldots+b_n-1}\int_{I^{n-1}}r^{-\e}\eta_{n-1}^{\e-b_n}\eta_{n-1}^{b_n-1}
\prod_{i=1}^{n-1}(1-\xi_i)^{b_i-1}\bigwedge_{i=1}^{n-1}d\xi_i\nah\\
&\leq& r^{a_n}(1-r)^{b_1+\ldots+b_n-1}\int_{I^{n-1}}r^{-\e}\eta_{n-1}^{\e-1}
\prod_{i=1}^{n-1}(1-\xi_i)^{b_i-1}\bigwedge_{i=1}^{n-1}d\xi_i\nah\\
&=&r^{a_n}(1-r)^{b_1+\ldots+b_n-1}r^{-\e}\prod_{i=1}^{n-1}\frac{\Gamma(\e)\Gamma(b_i)}{\Gamma(\e+b_i)}.\nah
\end{eqnarray}
As we did in section 2, this product of quotients of Gamma functions is meromorphic on the interval
$[0,b_n]$ except at zero where it has a pole of order $n-1$, and so there exists a constant $C$ such that
$$\prod_{i=1}^{n-1}\frac{\Gamma(\e)\Gamma(b_i)}{\Gamma(\e+b_i)}\leq C\e^{1-n}.$$
Hence
$$h(r)\ll r^{a_n}(1-r)^{b_1+\ldots+b_n-1}r^{-\e}\e^{1-n}.$$
Now
\begin{itemize}
\item[$\centerdot$] If $r\leq e^{-\frac{1}{b_n}}$, then $\ds\frac{1}{\ln\left(\frac{1}r\right)}\leq b_n$. In this case, we choose $\ds\e=\frac1{\ln(\frac1r)}$ and obtain
\begin{equation}\label{less}
h(r)\ll r^{a_n}(1-r)^{b_1+\ldots+b_n-1}e\left(\ln\left(\frac 1r\right)\right)^{n-1}.
\end{equation}
\item[$\centerdot$] If $r\geq e^{-\frac{1}{b_n}}$, we choose $\e=b_n$
and obtain
\begin{eqnarray}\label{big}
h(r)&\ll& r^{a_n}(1-r)^{b_1+\ldots+b_n-1}r^{-b_n}b_n^{1-n}\nah\\
&\leq& r^{a_n}(1-r)^{b_1+\ldots+b_n-1}eb_n^{1-n}.
\end{eqnarray}
\end{itemize}
Combining (\ref{less}) and (\ref{big}) together yields
\begin{eqnarray*}
h(r)&\ll& r^{a_n}(1-r)^{b_1+\ldots+b_n-1}\left(e\left(\ln\left(\frac 1r\right)\right)^{n-1}+eb_n^{1-n}\right)\\
&\ll& r^{a_n}(1-r)^{b_1+\ldots+b_n-1}\left(\ln\left(\frac 1r\right)\right)^{n-1},\textrm{ for all }0<r<1.
\end{eqnarray*}
\item[\bf{Case 4.}] Assume there exists $k$ such that $a_k>a_{k+1}=\ldots=a_n=a$. Let $F(r)$ be the convolution product of $f_1,f_2,\ldots,f_{k+1}$ and $G(r)$ be the convolution product of the rest
namely $f_{k+2}\ldots,f_n$. From the previous discussion it is clear that
$$F(r)\ll r^a(1-r)^{b_1+\ldots+b_{k+1}-1}$$
and
$$G(r)\ll r^a(1-r)^{b_{k+2}+\ldots+b_{n}-1}\left(\ln\left(\frac er\right)\right)^{n-k-2}.$$
Denote by $b=b_1+\ldots+b_{k+1}$, $d=b_{k+2}+\ldots+b_{n}$ and $n-k-1=l$. The case $l=1$ has been treated previously. So assume $l>1$. We see that
\begin{eqnarray*}
h(r)&=&(F*G)(r)\\
&\ll&\int_r^1(r/t)^a(1-r/t)^{b-1}t^a(1-t)^{d-1}\left(\ln\left(\frac et\right)\right)^{l-1}\frac{dt}{t}\nah\\
&\leq &r^a\int_r^1(t-r)^{b-1}t^{-b}(1-t)^{d-1}\left(\ln\left(\frac et\right)\right)^{l-1}\frac{dt}{t}.\nah
\end{eqnarray*}
Now the change of variables $t=u+r-ur$ leads to $t-r=u(1-r)$, $1-t=(1-u)(1-r)$, $dt=(1-r)du$
and
\begin{eqnarray*}
h(r)&\ll& r^a\int_0^1(1-r)^{b-1}u^{b-1}t^{-b}(1-r)^{d-1}(1-u)^{d-1}
\left(\ln\left(\frac{e}{t}\right)\right)^{l-1}(1-r)\,du\\
&\leq& r^a(1-r)^{b+d-1}\int_0^1u^{b-1}t^{-b}(1-u)^{d-1}\left(\ln\left(\frac et\right)\right)^{l-1}\,du.
\end{eqnarray*}
Noting that $t\geq u$ and $r >0$ and choosing an arbitrary $0<\e\leq b$ implies
\begin{eqnarray*}
h(r)&\ll& r^a(1-r)^{b+d-1}\int_0^1u^{b-1}t^{-\e}t^{\e-b}(1-u)^{d-1}\left(\ln\left(\frac et\right)\right)^{l-1}\,du\\
&\leq& r^a(1-r)^{b+d-1}\int_0^1u^{b-1}r^{-\e}u^{\e-b}(1-u)^{d-1}\left(\ln\left(\frac et\right)\right)^{l-1}\,du\\
&\leq&r^a(1-r)^{b+d-1}r^{-\e}\int_0^1u^{b-1}u^{\e-b}(1-u)^{d-1}\left(\ln\left(\frac eu\right)\right)^{l-1}\,du\\
&\leq&r^a(1-r)^{b+d-1}r^{-\e}\int_0^1u^{\e-1}(1-u)^{d-1}\left(\ln\left(\frac eu\right)\right)^{l-1}\,du.
\end{eqnarray*}
Let $\ds H_j(\e)=\int_0^1u^{\e-1}(1-u)^{d-1}(\ln u)^jdu$. This is the $j-$th order derivative
of the beta function $B(\e,d)$ as a function of $\e$, and $B(\e,d)$ is holomorphic on $(-1,\infty)$ except at zero where it has a simple pole with residue 1. This is easy to verify. So $\e^{j+1}H_j(\e)$ will be
holomorphic on the interval $(-1,\infty)$. Observing that
$$\int_0^1u^{\e-1}(1-u)^{d-1}\left(\ln\left(\frac eu\right)\right)^{l-1}\,du$$ is a linear sum of the derivatives of order less than or equal to $l-1$
of the Beta function, we find $$\e^l\int_0^1u^{\e-1}(1-u)^{d-1}\left(\ln\left(\frac eu\right)\right)^{l-1}du$$ is bounded by
a constant $C$ in the interval $[0,b]$. Thus
$$h(r)\ll r^a(1-r)^{b+d-1}r^{-\e}\e^{-l}.$$ Now arguing as in Case 3,
if $r\leq\ds e^{-\frac1b}$, we choose $\e=\frac{1}{\ln\left(\frac1r\right)}$
and get $$h(r)\ll r^a(1-r)^{b+d-1}e\left(\ln\left(\frac1r\right)\right)^l$$
and if $r>\ds e^{-\frac1b}$, we let $\e=b$, and have
$$h(r)\ll r^a(1-r)^{b+d-1}\frac e{b^l}.$$ Combining these two cases,
we obtain $$h(r)\ll  r^a(1-r)^{b+d-1}\left(\ln\left(\frac{e}{r}\right)\right)^l.$$ This totally proves the Lemma A.
\end{itemize}

\section{ Proof of Lemma B } We recall (\ref{h})
\begin{eqnarray*}
h(r)&=&r^{a_n}(1-r)^{b_1+\ldots+b_n-1}\int_{I^{n-1}}\eta_{n-1}^{b_n-1}
\prod_{i=1}^{n-2}\eta_i^{b_{i+1}}\\
& &\prod_{i=1}^{n-1}y_i^{a_i-a_{i+1}-b_{i+1}}
\prod_{i=1}^{n-1}(1-\xi_i)^{b_i-1}\bigwedge_{i=1}^{n-1}\,d\xi_i.
\end{eqnarray*} To make the differentiation easier, we introduce some notation.
Let
$$A=a_n,\  B=b_1+\ldots+b_n,\  \eta=(\eta_1,\ldots,\eta_{n-1}),\
\xi=(\xi_1,\ldots,\xi_{n-1}),\ y=(y_1,\ldots,y_{n-1})$$
$$\alpha_i=a_i-a_{i+1}-b_{i+1}\textrm{ for }1\leq i\leq n-1,\
\beta_i=b_{i+1}\textrm{ for }1\leq i\leq n-2,\ \beta_{n-1}=b_n-1,$$
$$\beta=(\beta_1,\ldots,\beta_{n-1}),\ G(\xi)=\prod_{i=1}^{n-1}(1-\xi_i)^{b_i-1},\
d\xi=\bigwedge_{i=1}^{n-1}d\xi_i,\ J=I^{n-1}.$$ With this notation and the multi-index
notation like for example $y^{\alpha}=y_1^{\alpha_1}\ldots y_{n-1}^{\alpha_{n-1}}$, (\ref{h}) can be written as
\begin{equation}\label{newh}
h(r)= r^A(1-r)^{B-1}\int_Jy^{\alpha}\eta^{\beta}G(\xi)\,d\xi.
\end{equation}

Clearly the function $\eta^\beta G(\xi_i)$ is summable $d\xi$, and each $y_i=\eta_i+r(1-\eta_i)$ satisfies
$0<r\leq y_i<1$ for $0<r<1$. So one can differentiate under the integral sign with respect to $r$. But before we do that let us introduce some more notation
$$g_1(r)=r^A,\ g_2(r)=(1-r)^{B-1},\ u_i=y_i^{\alpha_i}\textrm{ for }1\leq i\leq n-1.$$ Rewriting (\ref{newh}) as
\begin{equation}\label{hgg}
h(r)= \int_Jg_1g_2u_1\ldots u_{n-1}\eta^{\beta}G(\xi)\,d\xi.
\end{equation}
Now differentiating under the integral sign,
we obtain
\begin{equation}\label{Dh}
h^{(k)}(r)=\sum \int_Jg_1^{(l_1)}g_2^{(l_2)}u_1^{(j_1)}\ldots u_{n-1}^{(j_{n-1})}\eta^{\beta}G(\xi)\,d\xi
\end{equation}
 where the summation
is taken over all $(n+1)-$tuples of non-negative integers $(l_1,l_2,j_1,\ldots j_{n-1})$ such that
$k=l_1+l_2+j_1+\ldots+j_{n-1}$. Further it easy to check the following
\begin{eqnarray*}
u_i^{(j_i)}(r) & = & \alpha_i(\alpha_i-1)\ldots (\alpha_i-j_i+1)y_i^{\alpha_i-j_i}(1-\eta_i)^{j_i},\\
g_1^{(l_1)}(r)& = & A(A-1)\ldots (A-l_1+1)r^{A-l_1},\\
g_2^{(l_2)}(r) & = & (B-1)(B-2)\ldots(B-l_2)(-1)^{l_2}(1-r)^{B-l_2-1}.\\
\end{eqnarray*}
Since $y_i\geq r$ and $0\leq\eta_i\leq 1$, the equalities above imply
\begin{eqnarray}\label{u}
u_i^{(j_i)}(r) &\ll& y_i^{\alpha_i}r^{-j_i}\\
\label{g1} g_1^{(l_1)}(r)&\ll& g_1(r)r^{-l_1}\\
\label{g2} g_2^{(l_2)}(r)&\ll& (1-r)^{B-k-1}=g_2(r)(1-r)^{-k}
\end{eqnarray}
where the last inequality is obtained because  $0\leq l_2\leq k$. From (\ref{u}), (\ref{g1}) and (\ref{g2})  we deduce that
\begin{eqnarray*}
g_2^{(l_2)}(r)g_1^{(l_1)}(r)u_1^{(j_1)}(r)\ldots u_{n-1}^{(j_{n-1})}(r)&\ll& g_2(r)(1-r)^{-k}g_1(r)u_1(r)\ldots \\
&&\ldots u_{n-1}(r)r^{-l_1-j_1-\ldots-j_{n-1}}\\
&\ll& r^{-k}(1-r)^{-k}g_2(r)g_1(r)u_1(r)\ldots u_{n-1}(r).\end{eqnarray*}
Multiplying both sides by $\eta^{\beta}G(\xi)d\xi$ and integrating over $J$ yield
$$h^{(k)}(r)\ll r^{-k}(1-r)^{-k}h(r)$$ and by the Lemma A,
$$h(r)\ll r^A(1-r)^{B-1}\left(\ln\left(\frac er\right)\right)^{n-1}.$$ Hence we have
$$h^{(k)}(r)\ll r^{A-k}(1-r)^{B-k-1}\left(\ln\left(\frac{e}{r}\right)\right)^{n-1}.$$
This proves the Lemma B.


{\small
\begin{thebibliography}{9}
\bibitem{cr} \u Z. \u Cu\u ckovi\'c and N. V. Rao, Mellin transform,
monomial Symbols, and commuting Toeplitz operators, \emph{J.
Funct. Anal.} {\bf{154}} (1998), 195-214.
\bibitem {lz} I. Louhichi, L. Zakariasy, On Toeplitz operators with
  quasihomogeneous symbols, \emph{Arch. Math.} {\bf {85}}
   (2005), 248-257.
\bibitem{lsz} I. Louhichi, E. Strouse, L. Zakariasy, Products of
Toeplitz operators on the Bergman space, \emph{Integral Equations
Operator Theory} {\bf{54}} (2006), 525-539.
\bibitem{l} I. Louhichi, Powers and roots of Toeplitz operators,
\emph{Proc. Amer. Math. Soc.} {\bf{135}}, (2007), 1465-1475.
\bibitem{cl}  \u Z. \u Cu\u ckovi\'c and I. Louhichi, Finite rank commutators and
semicommutators of Toeplitz operators with quasihomogeneous
symbols. \emph{ Complex Analysis and Operator Theory.} Volume {\bf{2}}, Number 3 (2008), 429-439.
\bibitem{lr} I. Louhichi and N. V. Rao, Bicommutants of Toeplitz
operators, \emph{Arch. Math.} {\bf{91}} (2008), 256-264.
\bibitem{lry} I. Louhichi, N. V. Rao and A. Yousef, Two questions on
products of Toeplitz operators on the Bergman space, \emph {Complex Analysis and Operator Theory.}
Volume {\bf{3}}, Number 4 (2009), 881-889.


\end{thebibliography}}
\end{document}